\documentclass{article}

\usepackage{latexsym}
\usepackage{amsmath}
\usepackage{amssymb}

\newenvironment{itemise}{\begin{itemize}}{\end{itemize}}

\newcommand{\brak}[1]{\ensuremath{\left[#1\right]}}

\newcommand{\paren}[1]{\ensuremath{\left(#1\right)}}

\newcommand{\pr}[1]{\mathrm{Pr}\paren{#1}}

\newcommand{\goesto}{\rightarrow}

\newcommand{\reals}{\mathbb{R}}
\newcommand{\td}{\ensuremath{T_{\delta}}}
\newcommand{\tdd}{\ensuremath{T_{\delta 1}}}

\newcommand{\wrt}{\,\mathrm{d}}

\begin{document}

\title{A Note on Walking Versus Waiting}
\author{Anthony B. Morton}
\date{February 2008}
\maketitle

To what extent is a traveller (called Justin, say) better off to wait for a bus rather than just start walking---particularly when the bus headway is of a similar order of magnitude to the walking time and Justin does not know the precise arrival time of the bus?
The recent analysis by Chen et al \cite{chen2008} goes some way toward answering this question; however, there are additional valuable insights that can be gained from their approach.

In essence this is a decision-theoretic problem, where the aim is to compare various strategies for getting to one's destination.
Following Chen et al \cite{chen2008}, it is supposed that the comparison hinges solely on the expectation of overall travel time, and there is no additional penalty attached to walking.
We may list some of the available strategies as follows:
\begin{itemise}
\item
Strategy A: Wait indefinitely until a bus arrives.
\item
Strategy B: Begin walking immediately and do not wait for a bus.
\item
Strategy C: Wait for a predetermined time interval $T_W$, and if no bus has arrived in this time then walk.
\end{itemise}
(Naturally, strategies A and B can themselves be thought of as the extreme limiting cases of strategy C as $T_W \goesto \infty$ and $T_W \goesto 0$ respectively.
However, we distinguish them for the sake of discussion.)

Following the nomenclature of \cite{chen2008}, let $d$ denote the journey distance, $v_w$ Justin's average walking speed and $v_b$ the average speed of a bus.
The overall travel time is a random variable and will be denoted $T_t$.

As indicated in \cite{chen2008}, if the decision is purely between strategies A and B, the problem has an easy, one might say trivial, solution: if the expected arrival time for the bus is $T$ minutes in the future, Justin should wait if $T \leq d / v_w - d / v_b$ and walk otherwise.
The more interesting questions are (a) could there exist values $T_W$ for which strategy C is superior to either A or B; and (b) does the presence of intermediate stops affect the conclusion?
The purpose of this note is to show, by direct analysis for (almost) arbitrary arrival probabilities, that the answer to (a) is indeed no in most practical cases, but that the answer to (b) is not quite as clear-cut as the authors in \cite{chen2008} argue.

\section*{The Case of No Intermediate Stops}

The following expression is provided in \cite{chen2008} for Justin's expected travel time under strategy C, where $p(t)$ is the density function for the probability a bus arrives $t$ minutes from now:
\begin{equation}
E[T_t] = \int_0^{T_W} \paren{\frac{d}{v_b} + \tau} p(\tau) \wrt\tau
   + \paren{1 - \int_0^{T_W} p(\tau) \wrt\tau} \paren{\frac{d}{v_w} + T_W}.
\label{eq:exptt}
\end{equation}
In \cite{chen2008}, this expression is used to determine the waiting time $T_W$ such that $E[T_t]$ is equal to the time $d / v_w$ required under strategy B.
However, it is also of interest to optimise (\ref{eq:exptt}) directly.
The derivative of (\ref{eq:exptt}) with respect to $T_W$ is
\begin{equation}
\frac{\partial}{\partial T_W} E[T_t] = \paren{1 - \int_0^{T_W} p(\tau) \wrt\tau} - \paren{\frac{d}{v_w} - \frac{d}{v_b}} p(T_W).
\label{eq:expttd}
\end{equation}
The first term in parentheses is the complementary c.d.f.\ of the arrival time for the bus; that is, the probability that the bus has \emph{not} arrived by time $T_W$.
Denote this function by $R(t)$:
\begin{equation}
R(t) = 1 - \int_0^t p(\tau) \wrt\tau.
\label{eq:rdef}
\end{equation}
Also for brevity, define \td\ as the difference between the travel time on foot and the travel time by bus:
\begin{equation}
\td = \frac{d}{v_w} - \frac{d}{v_b}.
\label{eq:tddef}
\end{equation}
Then the derivative (\ref{eq:expttd}) may be expressed as
\begin{equation}
\frac{\partial}{\partial T_W} E[T_t] = R(T_W) - \td\, p(T_W).
\label{eq:expttd2}
\end{equation}
Zeros of (\ref{eq:expttd2}) occur at points where
\begin{equation}
\frac{1}{\td} = \frac{p(T_W)}{R(T_W)} = \lambda(T_W).
\label{eq:expttopt}
\end{equation}
In reliability theory where one is concerned with failure-time distributions, the function $\lambda(t) = p(t) / R(t)$ is known as the \emph{hazard rate}: the density for the probability of failure given a component has not failed up to time $t$ \cite{billinton1992}.
But $\lambda(t)$ also has a natural interpretation for arrival times, where it may be termed the \emph{appearance rate} for the distribution $p(t)$.

By definition $\lambda(t)$ is nonnegative, but depending on the underlying distribution it may behave in a variety of ways: it may increase, decrease, attain local maxima or minima, or be constant.
More precisely, for any $T \in \reals_+ \cup \{\infty\}$ any nonnegative real measurable function on $(0,T)$ can be the appearance rate for a possible arrival time distribution $p(t)$; if $T$ is finite then take $p(t) = R(t) = 0$ for $t \geq T$.
An increasing $\lambda(t)$ is characteristic of many waiting scenarios; it represents the case where the conditional probability increases with time spent waiting.
By contrast, a constant $\lambda(t)$ characterises the exponential arrival time distribution, where time spent waiting does not affect the conditional likelihood of a bus arriving.

Equation (6) indicates that $E[T_t]$ as a function of $T_W$ is stationary at any point where the appearance rate $\lambda(T_W)$ equals the reciprocal of \td.
To determine whether this is a minimum or a maximum, one calculates the second derivative (assuming for simplicity that $p(t)$ is dfferentiable)
\begin{equation}
\frac{\partial^2}{\partial T_W^2} E[T_t] = - p(T_W) - \td\, p'(T_W)
\label{eq:expttdd}
\end{equation}
and inspects the sign of this quantity.
Since $p(t)$ and \td\ are both positive (and $p(T_W)$ is nonzero by (\ref{eq:expttopt})) this is equivalent to the sign of
\begin{equation}
- \paren{\frac{1}{\td} + \frac{p'(T_W)}{p(T_W)}}
   = - \paren{\frac{p(T_W)}{R(T_W)} + \frac{p'(T_W)}{p(T_W)}}
   = - \frac{\lambda'(T_W)}{\lambda(T_W)} = - \td\, \lambda'(T_W).
\label{eq:optsign}
\end{equation}
Accordingly, any stationary point of $E[T_t]$ is a minimum when it coincides with a falling appearance rate, and a maximum when it coincides with a rising appearance rate.

It is thus established that, when faced with a typical probability distribution for bus arrivals having a rising appearance rate, the optimal strategy reduces to a choice between strategies A and B, and hence to a straightforward comparision of expected arrival time with \td\ as indicated above.

For example, if it is known that buses run punctually every $T$ minutes but one does not know the actual arrival times, the uncertainty in arrival time can be modelled as a uniform probability density $p(t)$ between $0$ and $T$, with mean $T / 2$.
The corresponding appearance rate is
\begin{equation}
\lambda(t) = \frac{1}{T - t}, \qquad t \in (0,T)
\label{eq:auniform}
\end{equation}
an increasing function of $t$.
One may then distinguish three cases:
\begin{enumerate}
\item
If $T < \td$, the expected travel time $E[T_t]$ has no stationary points, and is in fact a decreasing function of wait time.
In this case strategy A is optimal: one should wait for the next bus.
\item
If $\td < T < 2\td$, the expected travel time has one stationary point at $T_W^* = T - \td$, which is a maximum.
The expected travel time increases with waiting time up to $T_W^*$ and then decreases.
But since the mean bus arrival time is $T / 2 < \td$, the waiting strategy A is still preferable to an initial decision to walk.
\item
If $T > 2\td$, then expected travel time again reaches a maximum at $T_W^* = T - \td$.
But now, it is better to follow strategy B and not wait for a bus at all.
\end{enumerate}
Perhaps the most interesting situation here is the marginal case where $T = 2\td$: the interval between buses is twice the difference between walking time and bus travel time.
In this case it is equally good (or bad) to wait as to walk, but having decided to wait, one should not give up and walk, as this will result in a worse outcome (averaged over all eventualities).

Not all practically realisable probability distributions have rising appearance rates for all $t$.
Consider, for example, the situation where Justin has arrived at the bus stop one minute after the (known) scheduled arrival time for the bus, but knows that buses sometimes run up to five minutes late.
Under reasonable assumptions this leads to a falling $\lambda(t)$ for the first four minutes, and if $p(0)$ is sufficiently large and the service relatively infrequent, this will result in an optimal time $T_W^*$ for which Justin should remain at the bus stop and then start walking.
However, if the next bus can be counted on to arrive in a time comparable to \td, it will still be advantageous to keep waiting for the next bus, and the lazy mathematician still wins.

One may also contemplate the special case where $\lambda(t) \equiv \lambda$ is constant, and bus arrivals are a Poisson process.
The mean time between arrivals is $T = 1 / \lambda$; if this is greater than \td\ then it is better to walk, and if less then it is better to wait, as is intuitively clear.
More interesting is the case where $T$ happens to just equal \td: then, the derivative (\ref{eq:expttd2}) vanishes and the expected travel time is equal to the walking time $d / v_w$, \emph{regardless} of the time spent waiting.
What this means is, no matter how long Justin decides in advance to wait, over a large number of journeys he can expect to spend no more or less time than if he had decided to walk at the outset.
Nonetheless, Poisson arrivals are generally an unrealistic assumption except in some cases of very frequent service, where it is generally better to catch the bus in any event.

\section*{Intermediate Stops}

In \cite{chen2008} it is argued that the conclusion is unchanged when Justin has the additional option of walking to an intermediate stop and catching the bus there.
If, for example, one re-evaluates the decision between strategies A and B at a distance $d_1$ from the start of the journey, then the walking time is $(d - d_1) / v_w$ and the bus travel time is $(d - d_1) / v_b$.
However, the expected waiting time for the bus (assuming it does not pass by \emph{en route}) is now reduced by $T_1 = d_1 / v_w - d_1 / v_b$, the time that passes while walking less the additional time taken by the bus to reach the next stop.
The terms involving $d_1$ cancel, and one is left with a decision identical to the original one.
Under these circumstances, one would ``rather save energy'' and act as though the intermediate stop did not exist, since the outcome is the same.
This conclusion is also intuitively evident.

The authors argue that the same reasoning applies to walk-and-wait strategies analogous to strategy C, at least when the wait time $T_W$ at the next stop is chosen to make the expected travel time equal to the walking time.
However, there appears to be a circularity involved here---if the expected travel time is fixed as $d / v_w$ \emph{a priori}, then naturally it will be observed to be the same whether or not Justin walks to another stop.
Again, some more insight into the problem is gained by seeking stationary points of $E[T_t]$ with respect to the free variables.

Suppose strategy C is modified so that one walks a distance $d_1$ to another stop, then waits for a maximum time $T_W$.
One may now distinguish two cases, depending on whether or not a bus passes by \emph{en route}.
Denote this event by $M$ (`miss the bus') and its non-occurrence by $\bar{M}$; the relevant probabilities are
\begin{equation}
\pr{M} = \int_0^{T_1} p(\tau) \wrt\tau, \qquad
\pr{\bar{M}} = 1 - \pr{M} = R(T_1), \qquad
T_1 = \frac{d_1}{v_w} - \frac{d_1}{v_b}.
\label{eq:prm}
\end{equation}
If no bus passes by \emph{en route}, then the expected travel time is similar to (\ref{eq:exptt}), but with a time offset:
\begin{align}
E[T_t | \bar{M}] = \frac{d_1}{v_w}
  &+ \int_0^{T_W} \paren{\frac{d - d_1}{v_b} + \tau} \frac{p(\tau + T_1)}{\pr{\bar{M}}} \wrt\tau \nonumber \\
  &+ \paren{1 - \int_0^{T_W} \frac{p(\tau + T_1)}{\pr{\bar{M}}} \wrt\tau} \paren{\frac{d - d_1}{v_w} + T_W}.
\label{eq:expttintbm}
\end{align}
Even if the bus does pass by, it is possible that Justin may catch it anyway, if he is vigilant and there is a stop close enough, or if the bus is delayed, or if he takes a short cut.
Let the probability of catching a bus if one turns up on the way be $P_C$.
On the other hand, if Justin misses the bus it may be assumed there is no advantage in waiting for the next one (else he would have preferred strategy A at the outset).
Thus the expected travel time conditional on $M$ is
\begin{equation}
E[T_t | M] = P_C \int_0^{T_1} \paren{\tau + \frac{d}{v_b}} \frac{p(\tau)}{\pr{M}} \wrt\tau
   + (1 - P_C) \frac{d}{v_w}.
\label{eq:expttintm}
\end{equation}
(Note that within the integral, $\tau$ is the time at which the caught bus reaches the \emph{starting} point, rather than Justin's current position, since $p(\tau)$ is the p.d.f.\ of bus arrivals at a \emph{fixed} point on the route.
This is a subtle point and easily overlooked.
In terms of $\tau$, Justin's overall travel time if he catches the bus is $\tau$ plus the time taken by the bus to cover the entire journey, regardless of where Justin is when he catches it.)

The overall expected travel time is
\begin{equation}
E[T_t] = \pr{M} E[T_t | M] + \pr{\bar{M}} E[T_t | \bar{M}],
\label{eq:expttint}
\end{equation}
where $E[T_t | M]$ and $E[T_t | \bar{M}]$ are given by (\ref{eq:expttintm}) and (\ref{eq:expttintbm}) respectively.
Again, one can differentiate to find optimal values of wait time $T_W$ at the more distant stop.
(This is made easier by noting that $E[T_t | M]$ does not vary with $T_W$.)
One finds that
\begin{align}
\frac{\partial}{\partial T_W} E[T_t] &= R(T_W + T_1) - \tdd p(T_W + T_1)
\label{eq:expttd1} \\ \mbox{and} \qquad
\frac{\partial^2}{\partial T_W^2} E[T_t] &= - p(T_W + T_1) - \tdd p'(T_W + T_1),
\label{eq:expttdd1}
\end{align}
where
\begin{equation}
\tdd = (d - d_1) \paren{\frac{1}{v_w} - \frac{1}{v_b}} = \td - T_1
\label{eq:tdd}
\end{equation}
is the difference between travel time on foot and by bus starting from the more distant stop.
These are exactly equivalent to (\ref{eq:expttd2}) and (\ref{eq:expttdd}), with $T_W + T_1$ in place of $T_W$ and \tdd\ in place of \td.
Accordingly, the analysis of the previous section is essentially unchanged, as is the conclusion: that in most waiting scenarios the optimal value of $T_W$ is either zero or arbitrarily large.

But now there is an additional variable in the problem, because Justin will generally have a number of stops to choose from between the starting point and the destination.
For the sake of simplicity, suppose the stop spacing is small enough that $d_1$ can be regarded as a continuous variable.
Then we can try and optimise with respect to $d_1$, obtaining after a little work
\begin{equation}
\frac{\partial}{\partial d_1} E[T_t] = q^2 (d - d_1) \brak{(1 - P_C) p(T_1) - p(T_1 + T_W)}
\label{eq:expttd3}
\end{equation}
where $q = 1 / v_w - 1 / v_b$ is a constant.

The above discussion suggests directing attention to the cases $T_W = 0$ and $T_W \goesto \infty$.
If $T_W = 0$, then (\ref{eq:expttd3}) vanishes unless $P_C > 0$; that is, when there is a nonzero chance of catching a bus while walking to the next stop.
Assuming then that $P_C > 0$ and $T_W = 0$, (\ref{eq:expttd3}) is found to be zero when $p(T_1) = 0$ and negative when $p(T_1) > 0$.
(Recall that $T_1$ is proportional to $d_1$.)
The result is intuitively evident: as long as there is a chance of catching a bus on the way (however small), walking a distance $d_1$ to an intermediate stop always reduces the expected travel time $E[T_t]$ relative to strategy B, where one does not attempt to catch a bus.
The actual reduction can be quantified by integrating (\ref{eq:expttd3}) with respect to $d_1$.
Due to the factor $(d - d_1)$ in (\ref{eq:expttd3}), the marginal reduction in travel time is greatest at the start of the journey, all other things being equal.

Now consider the case $T_W \goesto \infty$.
One may presume that there is a limit to the time one must wait for a bus, so that $p(T_1 + T_W)$ becomes zero for $T_W$ sufficiently large.
It is also fair to assume that $P_C < 1$.
Then it follows that for $T_W$ sufficiently large, (\ref{eq:expttd3}) is zero when $p(T_1) = 0$, but is \emph{positive} when $p(T_1) > 0$.
It follows that walking to an intermediate stop only to wait indefinitely at that stop yields a worse expected outcome than waiting at the starting point, consistent with common sense.
Again, (\ref{eq:expttd3}) can be used to judge how much one is worse off due to the chance of missing the bus on the way.

The more interesting case to consider, of course, is a direct comparison between strategy A (waiting indefinitely at the starting point) and a strategy of walking but trying to catch a bus on the way.
From above, the best possible strategy of the latter variety is the case $T_W = 0$ and $d_1 = d$: a little manipulation of (\ref{eq:expttint}) yields the result
\begin{equation}
E_*[T_t] = \frac{d}{v_w} - P_C \int_0^{\td} (\td - \tau) p(\tau) \wrt\tau.
\label{eq:expttintopt}
\end{equation}
Under strategy A, of course, the expected travel time is just the bus travel time, plus the expected arrival time for the bus:
\begin{equation}
E_A[T_t] = \frac{d}{v_b} + \int_0^{\infty} \tau p(\tau) \wrt\tau.
\label{eq:exptta}
\end{equation}
The difference in expected travel time is then
\begin{equation}
E_A[T_t] - E_*[T_t] = P_C \td \paren{1 - R(\td)} + (1 - P_C) \int_0^{\td} \tau p(\tau) \wrt\tau
   + \int_{\td}^{\infty} \tau p(\tau) \wrt\tau - \td.
\label{eq:expttdiff}
\end{equation}
If this quantity is positive, it represents the expected time saving if one walks instead of waiting.
(It will be seen that all terms are in fact nonnegative apart from the last one, \td; which, however, is relatively large.)

For example, let $p(t)$ be the uniform distribution based on a bus headway $T$.
The mean arrival time is $T/2$ and hence $E_A[T_t] = d / v_b + T / 2$.
If $T > \td$, then
\begin{equation}
E_*[T_t] = \frac{d}{v_w} - P_C \frac{\td^2}{2T}
\label{eq:uexptt1}
\end{equation}
and this improves on $E_A[T_t]$ if
\begin{equation}
P_C > 2 \frac{T}{\td} - \paren{\frac{T}{\td}}^2.
\label{eq:upc1}
\end{equation}
As $T / \td$ increases from 1 to 2, the minimum probability $P_C$ decreases from 1 to zero; recall that when $T / \td > 2$ strategy A is worse than strategy B, which is equivalent to the above strategy with $P_C = 0$.

If $T < \td$, then
\begin{equation}
E_*[T_t] = \frac{d}{v_w} - P_C \paren{\td - \frac{T}{2}}
\label{eq:uexptt2}
\end{equation}
and the condition for this to improve on $E_A[T_t]$ reduces to $P_C > 1$, a contradiction.
Accordingly, strategy A is still the best choice for uniform arrival probabilities with headway $T < \td$.

\section*{Conclusion}

Closer analysis of the model and results in \cite{chen2008} for a traveller at a bus stop confirms that the lazy mathematician does indeed win in many cases.
The exact criteria for when it is better to walk can be formulated in terms of an appearance rate function, and if one disregards the existence of intermediate stops, an arrival probability with a rising appearance rate is sufficient to make waiting the best strategy.

It has also been shown, however, that allowing for the possibility of catching the bus at an intermediate stop can restore the walking strategy to optimality.
In the end, the decision whether to walk or wait is not always clear-cut, and relies on one's expectations of being able to catch a bus on the run, as well as of a bus turning up in the first place.

In the particular case where one knows the (uniform) headway but not the exact timetable, it has been clearly established that one should wait if the headway is less than the walking time (less bus travel time), and should walk if the headway is more than twice this much.
This leaves a substantial window where it may be better to wait or to walk, depending on one's confidence in being able to catch up to a passing bus.


\begin{thebibliography}{1}

\bibitem{chen2008}
J.G. Chen, S.D. Kominers, and R.W. Sinnott.
\newblock Walk versus wait: The lazy mathematician wins.
\newblock {\em arXiv.org Mathematics}, January 2008.
\newblock http://arxiv.org/abs/0801.0297.

\bibitem{billinton1992}
R.~Billinton and R.N. Allan.
\newblock {\em Reliability Evaluation of Engineering Systems}.
\newblock Plenum, second edition, 1992.

\end{thebibliography}
\end{document}